\documentclass[12pt]{amsart}

\usepackage{amssymb}
\usepackage{enumitem}
\usepackage[T1]{fontenc}

\makeatletter
\@namedef{subjclassname@2020}{%
  \textup{2020} Mathematics Subject Classification}
\makeatother

\newtheorem{thm}{Theorem}
\newtheorem{lem}{Lemma}
\newtheorem{prop}{Proposition}
\newtheorem{cor}{Corollary}

\theoremstyle{definition}
\newtheorem{defn}{Definition}
\newtheorem{rem}{Remark}

\newcommand{\N}{\mathbb{N}}
\newcommand{\R}{\mathbb{R}}
\newcommand{\C}{\mathbb{C}}
\newcommand{\T}{\mathbb{T}}
\newcommand{\D}{\mathbb{D}}

\newcommand{\Qbar}{\overline{Q}}
\newcommand{\Acal}{\mathcal{A}}

\begin{document}

\title[Nevanlinna--Pick norms]{Nevanlinna--Pick norms on scattered and Cantor spectra}

\author[P. Ohrysko]{Przemys\l aw Ohrysko}
\address{University of Warsaw\\
Faculty of Mathematics, Informatics, and Mechanics\\
\\Institute of Mathematics\\
Banacha 2, 02-097 Warsaw, Poland}
\email{p.ohrysko@gmail.com}

\author[M. Wojciechowski]{Micha\l{} Wojciechowski\textsuperscript{*}}
\thanks{\textsuperscript{*} This research was partially supported by the National Science Centre, Poland, and Austrian Science Foundation FWF joint CEUS programme. National Science Centre project no.~2020/02/Y/ST1/00072 and FWF project no.~I5231.}
\address{Institute of Mathematics\\
Polish Academy of Sciences\\
00-656 Warszawa, Poland}
\email{miwoj.impan@gmail.com}

\date{}

\begin{abstract}
We introduce Nevanlinna--Pick norms associated with finite families of characters on a commutative semisimple Banach algebra $A$ and study the class $NP_\infty$ of algebras for which all these norms coincide with the $\ell^\infty$ norm. Our positive result is topological: if $A\in NP_\infty$ and $K\subset \Delta(A)$ is compact scattered Hausdorff, then the restriction algebra $NP(A,K)$ is isometrically isomorphic to $C(K)$. In particular, if $\Delta(A)$ is compact scattered, then
\[
A\in NP_\infty
\iff
A\cong C(\Delta(A))
\quad\text{isometrically}.
\]
By contrast, we show that this rigidity fails in the presence of Cantor structure: for every compact Hausdorff space $S$ containing a Cantor subset, there exists a commutative semisimple unital Banach algebra $A\in NP_\infty$ with $\Delta(A)\cong S$ such that $A\neq C(S)$. As a consequence, for compact metrizable spectra, scatteredness is exactly the topological condition forcing rigidity.
\end{abstract}


\maketitle

\section{Introduction}

Let $A$ be a commutative semisimple Banach algebra (throughout the paper, we use standard facts from the theory of commutative Banach algebras without further comment; see Kaniuth~\cite{k} for background.) Given pairwise distinct characters
\[
\varphi_1,\dots,\varphi_n\in \Delta(A),
\]
the associated Nevanlinna--Pick norm on $\C^n$ is defined by
\[
\|(a_1,\dots,a_n)\|_{NP(A;\varphi_1,\dots,\varphi_n)}
=
\inf\bigl\{\|x\|_A:\widehat x(\varphi_i)=a_i\text{ for }i=1,\dots,n\bigr\}.
\]
We study the class $NP_\infty$ of algebras for which all these norms coincide with the supremum norm. Equivalently, every finite interpolation problem on the maximal ideal space has the minimal possible norm.

Our main purpose is to understand how this local interpolation property is reflected in the topology of the maximal ideal space. The answer is strikingly rigid on the scattered side and strikingly flexible in the presence of Cantor structure.

A classical topological precursor to this picture already appears in the theory of uniform algebras. Indeed, for a compact space $K$, Rudin proved that if $K$ is scattered, then the only uniform algebra on $K$ is $C(K)$, whereas if $K$ contains a copy of the Cantor set, then $C(K)$ contains a proper uniform algebra; see Proposition~3.6.4 in Dales and \"Ulger~\cite{du}. The point of the present paper is that an analogous scattered/Cantor principle persists far beyond the uniform-algebra setting. More precisely, for semisimple commutative Banach algebras $A\in NP_\infty$, compact scattered subsets of the maximal ideal space enforce local $C(K)$-rigidity, whereas Cantor structure permits nontrivial algebras in $NP_\infty$ that are not of the form $C(S)$.

The condition $A\in NP_\infty$ is stronger than the solvability of a single quotient interpolation problem: it requires the quotient norm associated with every finite family of distinct characters to agree exactly with the corresponding $\ell^\infty$ norm. Our aim is to determine how this simultaneous finite interpolation property constrains the topology of the maximal ideal space.

In particular, on compact metrizable spaces, where non-scatteredness is equivalent to containing a copy of the Cantor set, the positive and negative results fit into a common topological framework.

The first main result is local.

\medskip
\noindent\textbf{Local scattered rigidity.}
If $A\in NP_\infty$ and $K\subset \Delta(A)$ is compact scattered Hausdorff, then
\[
NP(A,K)\cong C(K)
\quad\text{isometrically}.
\]

Applying this with $K=\Delta(A)$ immediately yields the global theorem.

\medskip
\noindent\textbf{Global scattered rigidity.}
If $\Delta(A)$ is compact scattered, then
\[
A\in NP_\infty
\iff
A\cong C(\Delta(A))
\quad\text{isometrically}.
\]

Thus, on compact scattered spectra, the condition $A\in NP_\infty$ is exactly the condition that $A$ is the corresponding commutative $C^*$-algebra.

The second main result shows that this rigidity is genuinely topological. If $S$ contains a Cantor subset, then rigidity fails: there exists a commutative semisimple unital Banach algebra $A\in NP_\infty$ with $\Delta(A)\cong S$ but $A\neq C(S)$. The construction uses a concrete Sobolev-type model algebra on the square, its trace algebra on a product of fat Cantor sets, and a pullback over the ambient compact space.

The two arguments have separate chains of dependencies. Local rigidity follows from atomicity of measures on scattered compacta, the $NP_\infty$ interpolation property, and a Hahn--Banach density argument. The counterexample is built through the sequence
\[
\Acal\longrightarrow T\longrightarrow T'\longrightarrow A,
\]
where one passes from the Sobolev model algebra to a trace quotient, transports that quotient to a prescribed Cantor subset, and then forms the ambient pullback algebra.

Combining the positive and negative results, we obtain a clean picture for compact metrizable spectra: for such spaces, scatteredness is equivalent to countability, and it is exactly the topological condition that forces every algebra in $NP_\infty$ with that spectrum to be $C(K)$.

Our emphasis is on this topological perspective. We do not, however, know of any non-metrizable compact space $S$ that is non-scattered and contains no Cantor set, and for which it is known whether
\[
A \in NP_{\infty}
\qquad\text{and}\qquad
\Delta(A)=S
\]
imply that
\[
A=C(S).
\]

Some of the notions and results considered below have appeared in a somewhat different context in Dales and \"Ulger~\cite{du}. In particular, the $NP_{\infty}$ class of Banach algebras was defined in this book (see Remark \ref{rem:spc-equivalence}) under the name 'strongly pointwise contractive'. Moreover, a special case of Theorem \ref{thm:main-counterexample} i.e. the construction of Banach algebras $A\in NP_{\infty}$ such that $\Delta(A)\simeq [0,1]\text{ or }\Delta(A)\simeq\mathbb{T}$, but $A$ is not isomorphic to $C(K)$, is also present in the aforementioned monograph (see Remark \ref{rem:du-special-cases}).

\section{Nevanlinna--Pick norms}

Let $A$ be a commutative Banach algebra. We write $\Delta(A)$ for its Gelfand space, endowed with the weak$^*$ topology. Throughout the paper, we assume that $A$ is semisimple; that is, for every $0\neq x\in A$ there exists $\varphi\in\Delta(A)$ such that $\varphi(x)\neq 0$.

Each $x\in A$ defines its Gelfand transform $\widehat{x}:\Delta(A)\to\C$ by
\[
\widehat{x}(\varphi)=\varphi(x).
\]

\begin{lem}\label{lem:finite-character-surjectivity}
Let $\varphi_1,\dots,\varphi_n\in\Delta(A)$ be pairwise distinct. Then the linear map
\[
\Phi:A\to\C^n,
\qquad
\Phi(x)=\bigl(\varphi_1(x),\dots,\varphi_n(x)\bigr),
\]
is surjective.
\end{lem}

\begin{proof}
Distinct characters are linearly independent. Indeed, suppose that
\[
\sum_{j=1}^m c_j\varphi_j=0
\]
is a nontrivial relation with the smallest possible number of nonzero coefficients. After discarding zero coefficients, assume that every $c_j\neq 0$. A relation involving only one character is impossible because characters are nonzero, so $m\ge 2$. Since $\varphi_1\neq\varphi_m$, choose $a\in A$ such that $\varphi_1(a)\neq\varphi_m(a)$. Applying the relation to $ax$ and subtracting $\varphi_m(a)$ times the original relation evaluated at $x$ gives
\[
\sum_{j=1}^{m-1}c_j\bigl(\varphi_j(a)-\varphi_m(a)\bigr)\varphi_j(x)=0
\qquad (x\in A).
\]
This is a relation involving fewer characters, while the coefficient of $\varphi_1$ is nonzero, a contradiction. Hence $\varphi_1,\dots,\varphi_n$ are linearly independent. The coordinate functionals of $\Phi$ are therefore linearly independent, so $\operatorname{rank}\Phi=n$ and $\Phi$ is surjective.
\end{proof}

\begin{defn}\label{def:np-norm}
Let $A$ be a commutative semisimple Banach algebra and let
\[
\varphi_1,\dots,\varphi_n\in\Delta(A)
\]
be pairwise distinct. The associated Nevanlinna--Pick norm on $\C^n$ is defined by
\[
\|(a_1,\dots,a_n)\|_{NP(A;\varphi_1,\dots,\varphi_n)}
=
\inf\{\|x\|_A:\widehat{x}(\varphi_i)=a_i\text{ for }i=1,\dots,n\}.
\]

By Lemma~\ref{lem:finite-character-surjectivity}, the interpolation set is nonempty for every datum in $\C^n$; thus this is precisely the quotient norm induced by $\Phi$.

\end{defn}

\begin{rem}\label{rem:np-dominates-sup}
For every choice of pairwise distinct characters $\varphi_1,\dots,\varphi_n\in\Delta(A)$ and every $(a_1,\dots,a_n)\in\C^n$,
\[
\|(a_1,\dots,a_n)\|_{NP(A;\varphi_1,\dots,\varphi_n)}\ge \|(a_1,\dots,a_n)\|_\infty.
\]
Indeed, if $x\in A$ interpolates the values $a_i$, then
\[
\|x\|_A\ge \|\widehat{x}\|_\infty \ge \max_i |\widehat{x}(\varphi_i)|.
\]
\end{rem}

We say that $A$ belongs to the class $NP_\infty$ if all of its Nevanlinna--Pick norms are equal to the corresponding supremum norms.

\begin{rem}\label{rem:spc-equivalence}
Dales and \"Ulger call a Banach function algebra $A$ \emph{strongly pointwise contractive} if, for every $\varepsilon>0$, every finite family of distinct characters $\varphi_1,\dots,\varphi_n\in\Delta(A)$, and every $\zeta_1,\dots,\zeta_n\in\D$, there exists $f\in A$ with $\|f\|_A\le 1$ such that
\[
|\zeta_i-\widehat f(\varphi_i)|<\varepsilon
\qquad (i=1,\dots,n);
\]
see Dales and \"Ulger~\cite[Definition~3.5.6, pp.~200--201]{du}. This property is equivalent to $A\in NP_\infty$, since the closed unit ball for each finite quotient interpolation norm is the closure of the image of the closed unit ball of $A$. In the BSE terminology of~\cite{du}, where $C_{\mathrm{BSE}}(A)$ denotes the algebra of BSE functions on $\Delta(A)$, the same property is equivalent to
\[
\|\lambda\|_{A^*}=\|\lambda\|_1
\qquad
(\lambda\in L(A):=\operatorname{span}\Delta(A)),
\]
where $\|\sum_i\alpha_i\varphi_i\|_1=\sum_i|\alpha_i|$, and to the isometric identity
\[
\bigl(C_{\mathrm{BSE}}(A),\|\cdot\|_{\mathrm{BSE}}\bigr)
=
\bigl(C_b(\Delta(A)),\|\cdot\|_\infty\bigr);
\]
see Dales and \"Ulger~\cite[Theorem~5.4.5, p.~345]{du}.
\end{rem}

\begin{lem}\label{lem:complex-tietze}
Let $X$ be a normal topological space, let $E\subset X$ be closed, and let $f:E\to\C$ be continuous and bounded. Then $f$ has a continuous extension $F:X\to\C$ satisfying
\[
\|F\|_\infty=\|f\|_\infty.
\]
\end{lem}

\begin{proof}
Extend the real and imaginary parts of $f$ by the real-valued Tietze extension theorem, obtaining a continuous complex-valued extension $H$ of $f$. Put $R=\|f\|_\infty$. If $R=0$, take $F=0$. Otherwise compose $H$ with the radial retraction
\[
\rho_R(z)=
\begin{cases}
z,& |z|\le R,\\
Rz/|z|,& |z|>R.
\end{cases}
\]
Then $F:=\rho_R\circ H$ extends $f$ and has supremum norm $R$.
\end{proof}

\begin{thm}\label{thm:cstar}
Every commutative $C^*$-algebra belongs to $NP_\infty$.
\end{thm}

\begin{proof}
Let $A$ be a commutative $C^*$-algebra. By the commutative Gelfand--Naimark theorem, $A$ is isometrically isomorphic to $C_0(\Delta(A))$. Given distinct points $\varphi_1,\dots,\varphi_n\in\Delta(A)$ and values $a_1,\dots,a_n\in\C$, proceed as follows. If $A$ is unital, apply Lemma~\ref{lem:complex-tietze} to the closed subset $\{\varphi_1,\dots,\varphi_n\}$ of $\Delta(A)$. If $A$ is non-unital, apply the lemma to the one-point compactification $\Delta(A)\cup\{\infty\}$, assigning the value $0$ at $\infty$. In either case, we obtain an element of $A$ with norm $\|(a_1,\dots,a_n)\|_\infty$ and the required interpolation values. Therefore all Nevanlinna--Pick norms are equal to the supremum norm.
\end{proof}

\section{Scattered spectra and rigidity}

\begin{defn}\label{def:np-k}
Let $A$ be a commutative semisimple Banach algebra and let $K\subset \Delta(A)$ be compact. We denote by $NP(A,K)$ the algebra of restrictions of Gelfand transforms to $K$, equipped with the quotient norm
\[
\|f\|_{NP(A,K)}=
\inf\{\|x\|_A:\widehat{x}|_K=f\}.
\]
Equivalently,
\[
NP(A,K)=A/I_K,
\qquad
I_K:=\{x\in A:\widehat{x}|_K=0\}.
\]
\end{defn}

\begin{defn}\label{def:scattered-space}
A Hausdorff space $K$ is called \emph{scattered} if every nonempty subset of $K$ has an isolated point in its relative topology.
\end{defn}

\begin{lem}\label{lem:np-k-stable}
If $A\in NP_\infty$ and $K\subset\Delta(A)$ is compact, then $NP(A,K)\in NP_\infty$.
\end{lem}

\begin{proof}
Let
\[
q:A\to E:=A/I_K=NP(A,K)
\]
be the quotient map, and let $\psi_1,\dots,\psi_m\in\Delta(E)$ be pairwise distinct. Define
\[
\varphi_j:=\psi_j\circ q\in\Delta(A)
\qquad (j=1,\dots,m).
\]
The characters $\varphi_j$ are pairwise distinct because $q$ is surjective. For $(a_1,\dots,a_m)\in\C^m$, the definition of the quotient norm gives
\begin{align*}
\|(a_j)\|_{NP(E;\psi_1,\dots,\psi_m)}
&=\inf_{\psi_j(f)=a_j} \|f\|_E \\
&=\inf_{\varphi_j(x)=a_j} \|x\|_A \\
&=\|(a_j)\|_{NP(A;\varphi_1,\dots,\varphi_m)}\\
& =\|(a_j)\|_\infty.
\end{align*}
Hence $E\in NP_\infty$.
\end{proof}

\begin{lem}\label{lem:scattered-measures-atomic}
Let $K$ be a compact scattered Hausdorff space, and let $\mu\in M(K)$. Then there exist a finite or countable index set $J$, pairwise distinct points $(x_j)_{j\in J}\subset K$, and scalars $(\alpha_j)_{j\in J}$ such that
\[
\mu=\sum_{j\in J} \alpha_j \delta_{x_j},
\qquad
\sum_{j\in J} |\alpha_j|<\infty.
\]
In particular,
\[
\|\mu\|=\sum_{j\in J} |\alpha_j|.
\]
\end{lem}

\begin{proof}
This is the classical Pe\l czy\'nski--Semadeni theorem for compact scattered spaces; see Pe\l czy\'nski--Semadeni~\cite[Main Theorem]{ps} or Semadeni~\cite{sem}.
\end{proof}

\begin{thm}\label{thm:local-scattered}
Let $A$ be a commutative semisimple Banach algebra such that $A\in NP_\infty$. If $K\subset\Delta(A)$ is compact scattered Hausdorff, then the restriction of the Gelfand transform induces an isometric isomorphism
\[
NP(A,K)\cong C(K).
\]
\end{thm}

\begin{proof}
Set
\[
E:=NP(A,K).
\]
By Lemma~\ref{lem:np-k-stable}, the algebra $E$ also belongs to $NP_\infty$. Since its elements are functions on $K$, every point of $K$ acts on $E$ by evaluation, and hence $K\subset\Delta(E)$. Moreover, the identity map
\[
R:E\to C(K),\qquad R(u)=u,
\]
is contractive.

We first show that the closed unit ball $B_E$ is dense in $B_{C(K)}$ with respect to $\|\cdot\|_\infty$. Let $\mu\in M(K)=C(K)^*$. The required norming identity is immediate when $\mu=0$, so assume $\mu\neq 0$. By Lemma~\ref{lem:scattered-measures-atomic}, write
\[
\mu=\sum_{j\in J}\alpha_j\delta_{x_j},
\qquad
\|\mu\|=\sum_{j\in J}|\alpha_j|,
\]
where $J$ is finite or countable and the points $x_j$ are pairwise distinct. Choose an increasing sequence of finite sets $J_N\subset J$ whose union is $J$; if $J$ is finite, take $J_N=J$ for all sufficiently large $N$.

For $j\in J$, put
\[
\varepsilon_j=
\begin{cases}
\overline{\alpha_j}/|\alpha_j|,&\alpha_j\neq 0,\\
1,&\alpha_j=0.
\end{cases}
\]
Fix $N\in\N$ and $\eta>0$. Since $E\in NP_\infty$, the datum $(\varepsilon_j)_{j\in J_N}$ has $NP$-norm $1$. By the definition of the infimum, there exists $v_{N,\eta}\in E$ such that
\[
v_{N,\eta}(x_j)=\varepsilon_j
\quad (j\in J_N),
\qquad
\|v_{N,\eta}\|_E\le 1+\eta.
\]
Set
\[
u_{N,\eta}:=\frac{v_{N,\eta}}{1+\eta}\in B_E.
\]
The contractivity of $R$ gives $|u_{N,\eta}(x)|\le 1$ on $K$. Therefore
\begin{align*}
|\mu(u_{N,\eta})|
&\ge \frac{1}{1+\eta} \sum_{j\in J_N} |\alpha_j|
-\sum_{j\in J\setminus J_N} |\alpha_j|.
\end{align*}
Letting first $N\to\infty$ and then $\eta\downarrow 0$, we obtain
\[
\sup_{u\in B_E} |\mu(u)| \ge \|\mu\|.
\]
The reverse inequality follows from the contractivity of $R$; hence
\[
\sup_{u\in B_E} |\mu(u)| =\|\mu\|.
\]
Since $B_E$ is convex and balanced, the Hahn--Banach separation theorem implies that this norming identity forces sup-norm density in $B_{C(K)}$. Consequently,
\[
\overline{B_E}^{\|\cdot\|_\infty}=B_{C(K)}.
\]

We next prove that $R$ is surjective. Let $f\in B_{C(K)}$ and put $r_0=f$. Inductively, suppose that
\[
\|r_{n-1}\|_\infty \le 2^{-(n-1)}.
\]
Then $2^{n-1}r_{n-1}\in B_{C(K)}$. By the density just proved, choose $u_n\in B_E$ such that
\[
\|u_n-2^{n-1}r_{n-1}\|_\infty<\frac{1}{2},
\]
and define
\[
r_n:= r_{n-1}-2^{-(n-1)}u_n.
\]
It follows that
\[
\|r_n\|_\infty<2^{-n}.
\]
Thus the series
\[
v:=\sum_{n=1}^\infty 2^{-(n-1)}u_n
\]
converges in $E$, because
\[
\sum_{n=1}^\infty 2^{-(n-1)}\|u_n\|_E\le 2.
\]
Its partial sums converge uniformly to $f$, since the corresponding remainders are $r_n$. As $R$ is continuous, the $E$-limit $v$ satisfies $R(v)=f$. By homogeneity, $R(E)=C(K)$.

Finally, $R:E\to C(K)$ is a continuous bijection between Banach spaces. After identifying $E$ with $C(K)$ via $R$, the bounded inverse theorem yields a constant $C>0$ such that
\[
\|g\|_E\le C\|g\|_\infty
\qquad (g\in C(K)).
\]
Let $f\in B_{C(K)}$. By the density of $B_E$, choose $u_j\in B_E$ with
\[
\|u_j-f\|_\infty\to 0.
\]
Then
\[
\|f\|_E
\le \|u_j\|_E+\|f-u_j\|_E
\le 1+C\|f-u_j\|_\infty.
\]
Letting $j\to\infty$ gives $\|f\|_E\le 1$. Hence
\[
B_{C(K)}\subset B_E.
\]
The reverse inclusion follows from the contractivity of $R$, so the two unit balls coincide. Therefore
\[
\|u\|_E=\|u\|_\infty
\qquad (u\in E),
\]
and $R$ is an isometric algebra isomorphism.
\end{proof}

\begin{thm}\label{thm:scattered-main}
Let $A$ be a commutative semisimple Banach algebra such that $\Delta(A)$ is compact scattered Hausdorff. Then the following are equivalent:
\begin{enumerate}[label=\textup{(\roman*)}]
\item $A\in NP_\infty$;
\item the Gelfand transform identifies $A$ isometrically with $C(\Delta(A))$.
\end{enumerate}
\end{thm}

\begin{proof}
The implication \textup{(ii)}$\Rightarrow$\textup{(i)} follows from Theorem~\ref{thm:cstar}. Conversely, assume \textup{(i)}. Applying Theorem~\ref{thm:local-scattered} with $K=\Delta(A)$, we obtain
\[
NP(A,\Delta(A))\cong C(\Delta(A))
\]
isometrically. Since semisimplicity gives
\[
I_{\Delta(A)}=\{0\},
\]
we have $NP(A,\Delta(A))=A$, and the conclusion follows.
\end{proof}

\section{Cantor-set counterexamples}

In this section we show that the rigidity result above is genuinely topological. Once the maximal ideal space contains a Cantor subset, the condition $A\in NP_\infty$ no longer forces $A$ to be of the form $C(S)$. The construction proceeds in four stages: we first introduce a concrete Sobolev-type model algebra $\Acal$ on the square, then pass to its trace quotient on a fat Cantor product, then transport this quotient to a prescribed Cantor subset of an arbitrary compact space, and finally form a pullback algebra on the ambient compact space.

\subsection{The model algebra on the square}

Let $Q=(0,1)^2$ and $\Qbar=[0,1]^2$.

\begin{defn}\label{def:w11}
For an open set $\Omega\subset\R^d$, the Sobolev space $W^{1,1}(\Omega)$ consists of all functions $u\in L^1(\Omega)$ whose first-order distributional derivatives $\partial_1u,\dots,\partial_du$ belong to $L^1(\Omega)$. It is equipped with the norm
\[
\|u\|_{W^{1,1}(\Omega)}
=
\|u\|_{L^1(\Omega)}
+
\sum_{j=1}^d\|\partial_j u\|_{L^1(\Omega)}.
\]
We write $\nabla u=(\partial_1u,\dots,\partial_du)$ for the weak gradient.
\end{defn}

\begin{defn}\label{def:square-algebra}
Define
\[
\Acal:=W^{1,1}(Q)\cap C(\Qbar),
\qquad
\|u\|_{\Acal}:=\|u\|_{L^\infty(Q)}+\|\nabla u\|_{L^1(Q)}.
\]
Multiplication is pointwise.
\end{defn}

\begin{prop}\label{prop:square-banach}
$(\Acal,\|\cdot\|_{\Acal})$ is a commutative unital Banach algebra and
\[
\|uv\|_{\Acal}\le \|u\|_{\Acal}\,\|v\|_{\Acal}
\qquad (u,v\in\Acal).
\]
\end{prop}

\begin{proof}
For $u,v\in\Acal$,
\[
\|uv\|_\infty\le \|u\|_\infty\|v\|_\infty.
\]
Moreover, in the distributional sense,
\[
\nabla(uv)=u\nabla v+v\nabla u,
\]
so
\[
\|\nabla(uv)\|_1\le \|u\|_\infty\|\nabla v\|_1+\|v\|_\infty\|\nabla u\|_1.
\]
This proves submultiplicativity.

For completeness, let $(u_n)$ be a Cauchy sequence in $\Acal$. Since $|Q|=1$,
\[
\|u_n-u_m\|_{L^1(Q)}
\le \|u_n-u_m\|_{L^\infty(Q)}.
\]
Hence $(u_n)$ is Cauchy both in $C(\Qbar)$ and in $W^{1,1}(Q)$. Let $u$ be its uniform limit and $v$ its $W^{1,1}$-limit. The two limits agree almost everywhere; hence $u=v$ almost everywhere and $u\in\Acal$. The definition of the norm then shows that $u_n\to u$ in $\Acal$.
\end{proof}

\begin{lem}\label{lem:spec-acal}
For every $u\in\Acal$,
\[
r_{\Acal}(u)=\|u\|_\infty.
\]
\end{lem}

\begin{proof}
We have
\[
\|u^n\|_{\Acal}\ge \|u^n\|_\infty=\|u\|_\infty^n,
\]
so
\[
\liminf_{n\to\infty} \|u^n\|_{\Acal}^{1/n}\ge \|u\|_\infty.
\]
On the other hand,
\[
\nabla(u^n)=nu^{n-1}\nabla u,
\]
and therefore
\[
\|u^n\|_{\Acal}
\le \|u\|_\infty^n+n\|u\|_\infty^{n-1}\|\nabla u\|_1.
\]
Taking $n$-th roots and letting $n\to\infty$ gives the reverse inequality. The spectral radius formula yields the claim.
\end{proof}

\begin{lem}\label{lem:char-bound-acal}
Every character on $\Acal$ is continuous with respect to the uniform norm and satisfies
\[
|\varphi(u)|\le \|u\|_\infty
\qquad (u\in\Acal).
\]
\end{lem}

\begin{proof}
For a character $\varphi$, the value $\varphi(u)$ belongs to the spectrum of $u$. Hence, by Lemma~\ref{lem:spec-acal},
\[
|\varphi(u)|\le r_{\Acal}(u)=\|u\|_\infty.
\]
\end{proof}

\begin{lem}\label{lem:dense-acal}
$\Acal$ is uniformly dense in $C(\Qbar)$.
\end{lem}

\begin{proof}
Polynomials belong to $\Acal$ and separate points of $\Qbar$, so the Stone--Weierstrass theorem applies.
\end{proof}

\begin{thm}\label{thm:delta-acal}
The map $x\mapsto \delta_x$ is a homeomorphism from $\Qbar$ onto $\Delta(\Acal)$.
\end{thm}

\begin{proof}
Let $\varphi\in\Delta(\Acal)$. By Lemma~\ref{lem:char-bound-acal}, $\varphi$ is continuous with respect to the uniform norm on $\Acal$. By Lemma~\ref{lem:dense-acal}, it extends uniquely to a continuous linear functional
\[
\widetilde\varphi:C(\Qbar)\to\C.
\]
If $f,g\in C(\Qbar)$ and $u_n\to f$, $v_n\to g$ uniformly with $u_n,v_n\in\Acal$, then $u_nv_n\to fg$ uniformly; continuity of $\widetilde\varphi$ therefore gives
\[
\widetilde\varphi(fg)
=\lim_n\varphi(u_nv_n)
=\lim_n\varphi(u_n)\varphi(v_n)
=\widetilde\varphi(f)\widetilde\varphi(g).
\]
Thus the extension is multiplicative. By the classical Gelfand--Kolmogorov theorem, $\widetilde\varphi$ is evaluation at some point of $\Qbar$, hence so is $\varphi$. The map $x\mapsto \delta_x$ is therefore a continuous bijection from the compact space $\Qbar$ onto the Hausdorff space $\Delta(\Acal)$, and is thus a homeomorphism.
\end{proof}

\begin{thm}\label{thm:interp-acal}
Let $x_1,\dots,x_n\in Q$ and let $a_1,\dots,a_n\in\C$ satisfy $|a_i|\le 1$. For every $\varepsilon>0$ there exists $u\in\Acal$ such that
\[
u(x_i)=a_i\quad (i=1,\dots,n),
\qquad
\|u\|_{\Acal}\le 1+\varepsilon.
\]
\end{thm}

\begin{proof}
Choose radii $r_i>0$ such that the closed discs $\overline{B(x_i,r_i)}$ are pairwise disjoint, contained in $Q$, and satisfy
\[
\pi\sum_{i=1}^n r_i\le \varepsilon.
\]
Define
\[
\phi_i(x):=\max\bigl(1-|x-x_i|/r_i,0\bigr).
\]
Then $\phi_i\in\Acal$, $\phi_i(x_i)=1$, and the supports of the functions $\phi_i$ are pairwise disjoint. Moreover,
\[
\|\nabla\phi_i\|_1=\pi r_i.
\]
Set
\[
u:=\sum_{i=1}^n a_i\phi_i.
\]
Since the supports are disjoint, $\|u\|_\infty\le 1$ and $u(x_i)=a_i$ for all $i$. Also,
\[
\|\nabla u\|_1\le \sum_{i=1}^n \|\nabla\phi_i\|_1\le \varepsilon.
\]
Thus $\|u\|_{\Acal}\le 1+\varepsilon$.
\end{proof}

\subsection{A trace algebra on a fat Cantor set}

Let $K\subset \Qbar$ be compact. Define the closed ideal
\[
I_K:=\{u\in\Acal:u|_K\equiv 0\}
\]
and the quotient Banach algebra
\[
T:=\Acal/I_K.
\]
Identifying $T$ with the range $R_K(\Acal)\subset C(K)$ via restriction, the quotient norm is
\[
\|g\|_T=\inf\{\|v\|_{\Acal}:v|_K=g\}.
\]

\begin{lem}\label{lem:bumps-vanish}
If $x\in\Qbar\setminus K$, then there exists $u\in I_K$ such that $u(x)\neq 0$.
\end{lem}

\begin{proof}
Choose $r>0$ such that $\overline{B(x,r)}\cap K=\varnothing$ and let
\[
\phi(\cdot)=\max\bigl(1-|x-\cdot|/r,0\bigr).
\]
Then $\phi\in\Acal$, $\phi|_K=0$, and $\phi(x)=1$.
\end{proof}

\begin{thm}\label{thm:delta-trace}
For every compact $K\subset\Qbar$, the maximal ideal space of $T=\Acal/I_K$ is naturally homeomorphic to $K$.
\end{thm}

\begin{proof}
By the quotient-spectrum formula,
\[
\Delta(T)\cong \{\varphi\in\Delta(\Acal):\varphi(I_K)=0\}.
\]
Using Theorem~\ref{thm:delta-acal}, the characters of $\Acal$ are exactly the evaluations $\delta_x$. Now $\delta_x(I_K)=0$ means precisely that $u(x)=0$ for every $u\in I_K$. By Lemma~\ref{lem:bumps-vanish}, this happens exactly when $x\in K$. Thus every character of $T$ is evaluation at a unique point of $K$, and the resulting map
\[
K\to \Delta(T),\qquad x\mapsto \delta_x,
\]
is a homeomorphism.
\end{proof}

Choose a fat Cantor set $F\subset [1/4,3/4]$ of positive Lebesgue measure, obtained from a binary Cantor-type construction
\[
F=\bigcap_{m=0}^\infty F_m,
\]
where $F_m$ is the union of $2^m$ pairwise disjoint closed intervals and, at stage $m\ge 1$, exactly $2^{m-1}$ open gaps are removed. Put
\[
K_0:=F\times F\subset Q.
\]

We shall show that the trace algebra $T=R_{K_0}(\Acal)$ is a proper subalgebra of $C(K_0)$.

\begin{lem}\label{lem:ac-1d}
If $H\in W^{1,1}(0,1)$, then $H$ has an absolutely continuous representative and
\[
\int_a^b |H'(t)|\,dt\ge |H(b)-H(a)|
\qquad (0\le a<b\le 1).
\]
\end{lem}

\begin{proof}
This is standard; see \cite{leo,adf}.
\end{proof}

\begin{lem}\label{lem:no-w11-extension}
There exists $h\in C(F)$ such that every continuous extension $H\in C([0,1])$ of $h$ fails to belong to $W^{1,1}(0,1)$.
\end{lem}

\begin{proof}

Consider the binary construction
\[
F=\bigcap_{m=0}^\infty F_m
\]
fixed above; at stage $m$, exactly $2^{m-1}$ open gaps are removed.
For each $m\ge 1$, let
\[
\mathcal I_m=\{I_{m,1},\dots,I_{m,2^m}\}
\]
be the family of basic closed intervals forming $F_m$, ordered from left to right, and define the locally constant function
\[
\sigma_m:F\to\{-1,1\}
\]
by setting
\[
\sigma_m\equiv (-1)^\nu
\qquad\text{on }I_{m,\nu}\cap F
\quad (\nu=1,\dots,2^m).
\]
Choose a strictly increasing sequence of integers
\[
1\le m_1<m_2<\cdots
\]
so rapidly that
\[
\sum_{\ell>k}2^{-m_\ell/2}\le \frac{1}{4}\,2^{-m_k/2}
\qquad (k\ge 1).
\]
Now define
\[
h(x):=\sum_{k=1}^\infty 2^{-m_k/2}\sigma_{m_k}(x)
\qquad (x\in F).
\]
The series converges uniformly on $F$, so $h\in C(F)$.

Let $G_{m_k,j}=(a_{k,j},b_{k,j})$ be one of the $2^{m_k-1}$ open gaps removed at stage $m_k$. Its endpoints belong to $F$. For every $i<k$, the points $a_{k,j}$ and $b_{k,j}$ lie in the same basic interval of level $m_i$, so
\[
\sigma_{m_i}(a_{k,j})=\sigma_{m_i}(b_{k,j}).
\]
At level $m_k$, however, the endpoints lie in adjacent basic intervals, so
\[
|\sigma_{m_k}(b_{k,j})-\sigma_{m_k}(a_{k,j})|=2.
\]
Therefore
\begin{align*}
|h(b_{k,j})-h(a_{k,j})|
&\ge 2\,2^{-m_k/2}-2\sum_{\ell>k}2^{-m_\ell/2}\\
&\ge 2\,2^{-m_k/2}-\frac{1}{2}\,2^{-m_k/2}
=\frac{3}{2}\,2^{-m_k/2}.
\end{align*}
In particular,
\[
|h(b_{k,j})-h(a_{k,j})|\ge 2^{-m_k/2}.
\]

Now let $H\in C([0,1])$ be any continuous extension of $h$. If $H\in W^{1,1}(0,1)$, then by Lemma~\ref{lem:ac-1d}, on each gap $G_{m_k,j}$ one has
\[
\int_{a_{k,j}}^{b_{k,j}} |H'(t)|\,dt
\ge |H(b_{k,j})-H(a_{k,j})|
=|h(b_{k,j})-h(a_{k,j})|
\ge 2^{-m_k/2}.
\]
Summing over all $2^{m_k-1}$ gaps removed at stage $m_k$, we obtain
\[
\int_0^1 |H'(t)|\,dt
\ge 2^{m_k-1}2^{-m_k/2}=2^{m_k/2-1},
\]
which tends to $\infty$ as $k\to\infty$. This contradiction proves the lemma.
\end{proof}

\begin{lem}\label{lem:slices}
Let $f\in W^{1,1}(Q)$. There exists a representative $\widetilde f$ of $f$ such that, for almost every $x\in(0,1)$, the function
\[
y \mapsto \widetilde f(x,y)
\]
belongs to $W^{1,1}(0,1)$ and has weak derivative $y\mapsto \partial_2f(x,y)$. If, in addition, $f$ has a continuous representative on $\Qbar$, then for almost every such $x$ the absolutely continuous representative of the Sobolev slice agrees everywhere with the continuous slice $y\mapsto f(x,y)$.
\end{lem}

\begin{proof}
The first assertion is the standard slicing property for Sobolev functions, obtained from Fubini's theorem and approximation by smooth functions in $W^{1,1}(Q)$; see \cite{leo,adf}. For the second assertion, the Sobolev slice and the continuous slice agree almost everywhere. Since both representatives are continuous, they agree everywhere.
\end{proof}

\begin{thm}\label{thm:trace-not-onto}
For $K_0=F\times F$ as above,
\[
R_{K_0}(\Acal)\subsetneq C(K_0).
\]
\end{thm}

\begin{proof}
Let $h\in C(F)$ be as in Lemma~\ref{lem:no-w11-extension}, and define
\[
g(x,y):=h(y)
\qquad ((x,y)\in K_0).
\]
Suppose that $g=f|_{K_0}$ for some $f\in\Acal$. By Lemma~\ref{lem:slices}, for almost every $x\in(0,1)$ the function $y\mapsto f(x,y)$ belongs to $W^{1,1}(0,1)$. Since $F$ has positive measure, we may choose $x_0\in F$ for which this holds. Then
\[
H(y):=f(x_0,y)
\]
is a continuous $W^{1,1}$ extension of $h$, contradicting Lemma~\ref{lem:no-w11-extension}. Thus $g\notin R_{K_0}(\Acal)$.
\end{proof}

\subsection{Transport and pullback}

Let $S$ be a compact Hausdorff space and suppose that $S$ contains a Cantor subset $K'$. Since $K_0$ is homeomorphic to the Cantor set, fix a homeomorphism
\[
\theta:K_0\to K'.
\]
Let $T=\Acal/I_{K_0}$ and define
\[
T':=\{g\circ\theta^{-1}:g\in T\}\subset C(K')
\]
with the transported norm.

Now set
\[
A:=\{F\in C(S):F|_{K'}\in T'\},
\qquad
\|F\|_A:=\max\bigl\{\|F\|_{C(S)},\|F|_{K'}\|_{T'}\bigr\}.
\]

\begin{prop}\label{prop:A-on-S}
$(A,\|\cdot\|_A)$ is a commutative unital Banach algebra and
\[
\|F\|_{C(S)}\le \|F\|_A.
\]
\end{prop}

\begin{proof}
The constant function $1$ belongs to $A$, since its restriction to $K'$ is the unit of $T'$. If $F,G\in A$, then $(FG)|_{K'}=(F|_{K'})(G|_{K'})\in T'$, so $FG\in A$ and
\[
\|FG\|_A\le \|F\|_A\,\|G\|_A.
\]
Let $(F_n)$ be a Cauchy sequence in $(A,\|\cdot\|_A)$. Then $(F_n)$ is Cauchy in $C(S)$ and $(F_n|_{K'})$ is Cauchy in $T'$. Hence there exist
\[
F\in C(S),\qquad w\in T'
\]
such that
\[
F_n\to F \text{ in } C(S),
\qquad
F_n|_{K'}\to w \text{ in } T'.
\]
Since the norm on $T'$ dominates the supremum norm on $K'$, the second convergence implies uniform convergence on $K'$, and therefore $w=F|_{K'}$. Thus $F\in A$ and $A$ is complete.
\end{proof}

Let
\[
\rho:A\to T',
\qquad
\rho(F)=F|_{K'},
\]
and put
\[
J:=\ker\rho=\{F\in A:F|_{K'}=0\}.
\]

\begin{lem}\label{lem:J-identification}
\begin{enumerate}[label=\textup{(\alph*)}]
\item The map induced by $\rho$ identifies $A/J$ isometrically with $T'$ as Banach algebras.
\item $J$ is naturally isometric to $C_0(S\setminus K')$.
\end{enumerate}
\end{lem}

\begin{proof}
For \textup{(a)}, let $w\in T'$. By Lemma~\ref{lem:complex-tietze}, $w$ has an extension $F\in C(S)$ with
\[
\|F\|_{C(S)}=\|w\|_\infty\le \|w\|_{T'}.
\]
Since $F|_{K'}=w\in T'$, we have $F\in A$, $\rho(F)=w$, and
\[
\|F\|_A
=
\max\{\|F\|_{C(S)},\|w\|_{T'}\}
=
\|w\|_{T'}.
\]
Thus $\rho$ is surjective. Moreover, every lift $G\in A$ of $w$ satisfies $\|G\|_A\ge \|w\|_{T'}$, while the lift just constructed has norm exactly $\|w\|_{T'}$. Hence the induced quotient map is isometric.

For \textup{(b)}, note that $J$ consists exactly of the continuous functions on $S$ vanishing on $K'$. Restriction therefore identifies $J$ with $C_0(S\setminus K')$, and for $F\in J$ one has
\[
\|F\|_A=\|F\|_{C(S)},
\]
so the identification is isometric.
\end{proof}

\begin{lem}\label{lem:c0-spectrum}
If $X$ is locally compact Hausdorff, then $\Delta(C_0(X))\cong X$ via point evaluations.
\end{lem}

\begin{proof}
Let $X^*=X\cup\{\infty\}$ be the one-point compactification of $X$. The unitization of $C_0(X)$ is naturally identified with $C(X^*)$. Hence every character on $C_0(X)$ extends uniquely to a character on $C(X^*)$ and is therefore evaluation at a point of $X^*$ by the Gelfand--Kolmogorov theorem. Evaluation at $\infty$ on $C(X^*)$ restricts to the zero functional on $C_0(X)$, whereas a character is nonzero. Thus the point $\infty$ cannot occur.
\end{proof}

\begin{thm}\label{thm:delta-A}
Every character on $A$ is a point evaluation at some point of $S$. Hence
\[
\Delta(A)\cong S.
\]
\end{thm}

\begin{proof}
Let $\chi\in\Delta(A)$.

If $\chi$ annihilates $J$, then it factors through the quotient
\[
A/J\cong T'.
\]
Since $T'$ is isometrically isomorphic to $T=\Acal/I_{K_0}$ via transport by $\theta$, Theorem~\ref{thm:delta-trace} identifies $\Delta(T')$ with $K'$. Hence $\chi$ is evaluation at a point of $K'$.

Assume now that $\chi$ does not annihilate $J$. Then the restriction
\[
\chi|_J:J\to\C
\]
is a nonzero character on $J\cong C_0(S\setminus K')$. By Lemma~\ref{lem:c0-spectrum}, there exists $s\in S\setminus K'$ such that
\[
\chi(j)=j(s)
\qquad (j\in J).
\]
By Lemma~\ref{lem:J-identification}\textup{(b)}, choose $j_0\in J$ with $j_0(s)=1$. Then for every $F\in A$ we have $Fj_0\in J$, so
\[
\chi(F)=\chi(F)\chi(j_0)=\chi(Fj_0)=(Fj_0)(s)=F(s).
\]
Thus $\chi$ is evaluation at $s$.

Thus every character of $A$ is a point evaluation on $S$. It remains to verify that distinct points induce distinct evaluations, or equivalently, that $A$ separates points of $S$.

Two distinct points of $K'$ are separated by $T'$ and hence by a lift through the surjective map $\rho$. If $s\in S\setminus K'$ and $t\in K'$, then, using $J\cong C_0(S\setminus K')$, choose $j\in J$ with $j(s)\neq 0$; since $j(t)=0$, the points are separated. Finally, two distinct points of $S\setminus K'$ are separated by $C_0(S\setminus K')\cong J$. Thus the map
\[
S\to \Delta(A),\qquad s\mapsto \delta_s,
\]
is a continuous bijection from the compact space $S$ onto the Hausdorff space $\Delta(A)$, and hence a homeomorphism.

\end{proof}

\begin{thm}\label{thm:interp-A}
Let $M=\{x_1,\dots,x_n\}\subset S$ and let $a_1,\dots,a_n\in\C$ satisfy $|a_i|\le 1$. For every $\varepsilon>0$ there exists $F\in A$ such that
\[
F(x_i)=a_i\quad (i=1,\dots,n),
\qquad
\|F\|_A\le 1+\varepsilon.
\]
\end{thm}

\begin{proof}
Split $M=M_1\cup M_2$ with $M_1=M\cap K'$ and $M_2=M\setminus K'$.

If $M_1=\varnothing$, define $\phi$ on the closed set $E:=K'\cup M_2$ by
\[
\phi|_{K'}=0,
\qquad
\phi(x_i)=a_i\quad (x_i\in M_2).
\]
The function $\phi$ is continuous on $E$, because the finite set $M_2$ is disjoint from the closed set $K'$ and hence every point of $M_2$ is isolated in $E$. Moreover, $\|\phi\|_{C(E)}\le 1$. By Lemma~\ref{lem:complex-tietze}, $\phi$ extends to $F\in C(S)$ with $\|F\|_{C(S)}\le 1$. Since $F|_{K'}=0\in T'$, we have $F\in A$ and $\|F\|_A\le 1$.

Assume now that $M_1=\{y_1,\dots,y_m\}$ with $m\ge 1$, and write $b_j:=a_i$ whenever $y_j=x_i$. Let $z_j=\theta^{-1}(y_j)\in K_0$. By Theorem~\ref{thm:interp-acal}, there exists $u\in\Acal$ such that
\[
u(z_j)=b_j\quad (j=1,\dots,m),
\qquad
\|u\|_{\Acal}\le 1+\varepsilon.
\]
Let $g:=u|_{K_0}\in T$ and $w:=g\circ\theta^{-1}\in T'$. Then
\[
\|w\|_{T'}=\|g\|_T\le \|u\|_{\Acal}\le 1+\varepsilon.
\]
Define $\phi$ on the closed set $E:=K'\cup M_2$ by
\[
\phi|_{K'}=w,
\qquad
\phi(x_i)=a_i\quad (x_i\in M_2).
\]
Again, $\phi$ is continuous on $E$, because every point of the finite set $M_2$ is isolated in $E$. Moreover, $\|\phi\|_{C(E)}\le 1+\varepsilon$. By Lemma~\ref{lem:complex-tietze}, there exists $F\in C(S)$ extending $\phi$ with $\|F\|_{C(S)}\le 1+\varepsilon$. Since $F|_{K'}=w\in T'$, we have $F\in A$ and $\|F\|_A\le 1+\varepsilon$.
\end{proof}

\begin{thm}\label{thm:proper-subalgebra}
$A\neq C(S)$.
\end{thm}

\begin{proof}
Since $T\neq C(K_0)$ by Theorem~\ref{thm:trace-not-onto}, also $T'\neq C(K')$. Choose $\psi\in C(K')\setminus T'$ and extend it to $\Psi\in C(S)$ by Lemma~\ref{lem:complex-tietze}. Then $\Psi\notin A$.
\end{proof}

\begin{thm}\label{thm:main-counterexample}
Let $S$ be a compact Hausdorff space containing a Cantor subset. Then there exists a commutative semisimple unital Banach algebra $A\in NP_\infty$ such that:
\begin{enumerate}[label=\textup{(\arabic*)}]
\item $\Delta(A)\cong S$;
\item $A\neq C(S)$.
\end{enumerate}
\end{thm}

\begin{proof}
Theorem~\ref{thm:delta-A} gives $\Delta(A)\cong S$, while Theorem~\ref{thm:interp-A} implies that $A\in NP_\infty$. Finally, Theorem~\ref{thm:proper-subalgebra} gives $A\neq C(S)$.
\end{proof}

\begin{rem}\label{rem:du-special-cases}
The special cases $S=[0,1]$ and $S=\T$ were established by Dales and \"Ulger~\cite[Examples~3.5.13 and~3.5.14, pp.~203--208]{du}. More precisely, Example~3.5.13 (pp.~203--205) constructs a proper natural unital Banach function algebra on $[0,1]$ that is strongly pointwise contractive, while Example~3.5.14 (pp.~205--208) gives such an algebra on $\T$. In the BSE language, these examples are revisited in~\cite[Examples~5.4.13 and~5.4.16, pp.~348--349]{du}. By Remark~\ref{rem:spc-equivalence}, both algebras belong to $NP_\infty$.
\end{rem}

By the Cantor--Bendixson theorem, a compact metrizable space is scattered if and only if it is countable. Moreover, every non-scattered compact metrizable space contains a homeomorphic copy of the Cantor set; see Kuratowski~\cite[\S~23.V and \S~36.V]{kur}. Together with Theorems~\ref{thm:scattered-main} and~\ref{thm:main-counterexample}, these facts yield the following corollary.
\begin{cor}\label{cor:metrizable-rigidity}
Let $K$ be a compact metrizable Hausdorff space. Then the following are equivalent:
\begin{enumerate}[label=\textup{(\roman*)}]
\item $K$ is countable;
\item $K$ is scattered;
\item for every semisimple commutative unital Banach algebra $A$ such that
\[
A\in NP_\infty
\qquad\text{and}\qquad
\Delta(A)=K,
\]
the Gelfand transform identifies $A$ isometrically with $C(K)$.
\end{enumerate}
\end{cor}

\section*{Declaration on the use of AI-assisted tools}
During the preparation of this revised manuscript, AI-assisted tools were used for English-language editing and to improve editorial consistency. All mathematical arguments and conclusions were checked by the authors, who take full responsibility for the content of the manuscript.

\end{document}